\documentclass[12pt,twoside]{amsart}
\usepackage{amsthm,amsfonts,amssymb,amscd,euscript}

\usepackage[matrix,arrow]{xy}

\author{Florin Ambro} 
\address{DPMMS, CMS\\
University of Cambridge,
Wilberforce Road, Cambridge CB3 0WB, UK.}
\email{f.ambro@dpmms.cam.ac.uk}

\newcommand{\isoto}{{\overset{\sim}{\rightarrow}}}

\newcommand{\C}{{\mathbb C}}
\newcommand{\G}{{\mathbb G}}
\newcommand{\Q}{{\mathbb Q}}
\newcommand{\Z}{{\mathbb Z}}

\newcommand{\calF}{{\mathcal F}}

\newcommand{\calL}{{\mathcal L}}

\newcommand{\calO}{{\mathcal O}}

\newcommand{\bA}{{\mathbf A}}
\newcommand{\bB}{{\mathbf B}}

\newcommand{\bM}{{\mathbf M}}

\newcommand{\fa}{{\mathfrak a}}

\newcommand{\Alb}{\operatorname{Alb}}
\newcommand{\Aut}{\operatorname{Aut}}

\newcommand{\codim}{\operatorname{codim}}

\newcommand{\Ext}{\operatorname{Ext}}

\newcommand{\Gal}{\operatorname{Gal}}
\newcommand{\Hom}{\operatorname{Hom}}
\newcommand{\calHom}{\operatorname{\mathcal Hom}}

\newcommand{\Isom}{\operatorname{Isom}}
\newcommand{\Ker}{\operatorname{Ker}}

\newcommand{\mult}{\operatorname{mult}}

\newcommand{\rank}{\operatorname{rank}}
\newcommand{\Sing}{\operatorname{Sing}}
\newcommand{\Spec}{\operatorname{Spec}}
\newcommand{\Supp}{\operatorname{Supp}}
\newcommand{\trace}{\operatorname{trace}}

\theoremstyle{plain}
\newtheorem{thm}{Theorem}[section]

\newtheorem{lem}[thm]{Lemma}

\newtheorem{prop}[thm]{Proposition}

\theoremstyle{definition}
\newtheorem{defn}[thm]{Definition}

\newtheorem{exmp}[thm]{Example}

\newtheorem{ack}{Acknowledgments}   

\theoremstyle{remark}

\begin{document}

\bibliographystyle{amsalpha+}
\title[The moduli b-divisor]
{The moduli b-divisor of an lc-trivial fibration}
%\date{April 8, 2003}

\begin{abstract} We study positivity properties of
the moduli (b-)divisor associated to a relative log
pair $(X,B)/Y$ with relatively trivial log canonical
class. 
\end{abstract}

\maketitle

%\tableofcontents

%%%%%%%%%%%%%%%%%%%%%%%%%%%%%%%%%%%%%%%%%%%%%%%%%%%%%%%
\setcounter{section}{-1}
%%%%%%%%%%%%%%%%%%%%%%%%%%%%%%%%%%%%%%%%%%%%%%%%%%%%%%%
%%% Document name: modbdiv.tex
%%% Last modified: August 13, 2003
%%%%%%%%%%%%%%%%%%%%%%%%%%%%%%%%%%%%%%%%%%%%%%%%%%%%%%%

%%%%%%%%%%%%%%%%%%%%%%%%%%%%%%%%%%
%%%%%%%%%%%%%%%%%%%%%%%%%%%%%%%%%%

\section{Introduction}

%%%%%%%%%%%%%%%%%%%%%%%%%%%%%%%%%%%
%%%%%%%%%%%%%%%%%%%%%%%%%%%%%%%%%%%

\footnotetext[1]{1991 Mathematics Subject Classification. 
Primary: 14J10, 14N30. Secondary: 14E30.}

In this paper we continue the study of {\em lc-trivial 
fibrations} $(X,B)/Y$, that is, roughly, relative
log pairs with relatively trivial log canonical class 
$K+B$ (see ~\cite{BP}). These type of fibrations are
expected to play a key role in inductive arguments in 
the (Log) Minimal Model Program.

Given an lc-trivial fibration $f\colon (X,B)\to Y$, 
there exists a canonical decomposition of Kodaira type
$$
K+B\sim_\Q f^*(K_Y+B_Y+M_Y),
$$
where $B_Y$ and $M_Y$ are $\Q$-divisors on $Y$, called
the {\em discriminant} and {\em moduli} $\Q$-divisors
(Kawamata~\cite{Kaw97,sba}).
The discriminant measures the singularities of the
log pair $(X,B)$ over codimension one points of $Y$,
whereas the moduli $\Q$-divisor is expected to define 
the rational map from $Y$ to the moduli space of the 
generic fibre. As explained in \cite{BP}, the following
properties are desirable for applications: Inversion of 
Adjunction and (Effective) Semi-ampleness. 
Inversion of Adjunction (or, equivalently, Shokurov's BP
Conjecture~\cite{Shokurov03}) was established in \cite{BP}: 
the log pairs $(X,B)$ and $(Y,B_Y)$ have the same type
of singularities if $Y$ is sufficiently high in its 
birational class. As for the moduli part, it is known 
that $M_Y$ is numerically effective (nef) if $Y$ is 
sufficiently high in its birational class 
(Kawamata~\cite{sba}). Semi-ampleness predicts that 
in fact the linear system $|kM_Y|$ is free of base points 
if $k$ is large and divisible and $Y$ is sufficiently high 
in its birational class. 

The main results of this paper are two partial answers 
to the semiampleness of the moduli part: 
(a) if $M_Y$ is numerically trivial on high enough models $Y$, 
then $M_Y$ is a torsion $\Q$-divisor (Theorem~\ref{torsion}); 
(b) if the horizontal part of $B$ is effective and $Y$ is 
sufficiently high, then there exists a contraction
$h\colon Y\to Z$ to a projective variety $Z$ and
a nef and big $\Q$-divisor $A$ on $Z$ such that 
$M_Y\sim_\Q h^*A$ (Theorem~\ref{klglobal}). 

One application of Theorem~\ref{klglobal} is a
logarithmic version of a result of Kawamata~\cite{minit}
(Abundance and Ueno's K Conjecture for minimal models
with numerically trivial canonical class):

\begin{thm}\label{ak} 
Let $(X,B)$ be a projective log variety with Kawamata 
log terminal singularities such that $K+B\equiv 0$. 
Then 
\begin{itemize}
\item[(1)] There exists a positive integer
$b$ such that $b(K+B)\sim 0$.
\item[(2)] The Albanese map $X\to \Alb(X)$
is a surjective morphism, with connected fibers.
Furthermore, there exist an \'etale covering 
$A'\to \Alb(X)$, a projective log variety 
$(F,B_F)$ and an isomorphism 
\begin{displaymath}
\xymatrix{
(X,B)\times_{\Alb(X)}A' \ar[rr]^\simeq\ar[dr] &  & 
(F,B_F)\times A' \ar[dl] \\
&   A'  & } 
\end{displaymath}
\end{itemize}
\end{thm}
By Theorem~\ref{klglobal} and Theorem~\ref{ak}.(1),
we can generalize our main result in \cite{NM}: modulo
the Log Minimal Model Program and Log Abundance for
smaller dimensional varieties, the Abundance Conjecture 
is reduced to the case of log minimal models of maximal 
nef dimension (Theorem~\ref{ab}).
Another application of Theorem~\ref{klglobal} is a
generalization of a result of Nakayama~\cite{Nakayama}
(see also Fujita~\cite{fuj3}):

\begin{thm} Let $(X,B)$ be a projective log variety 
with Kawamata log terminal singularities, let
$f\colon X\to Y$ be a contraction to a proper
normal variety $Y$ and let $\omega$ be a $\Q$-Cartier 
divisor on $Y$ such that
$$
K+B\sim_\Q f^*\omega.
$$
Then there exists a $\Q$-Weil divisor $B_Y$ such 
that $(Y,B_Y)$ is a log variety with Kawamata log
terminal singularities and 
$
\omega \sim_\Q K_Y+B_Y.
$
\end{thm}
The techniques we use are due to 
Fujita~\cite{fujita, fuj2}, Viehweg~\cite{wp, wpII} 
and Kawamata~\cite{koddim, minit}. 
In fact, Theorem~\ref{klglobal} can be deduced 
from \cite{minit}, Theorem 1.1 in case the generic fibre 
$X_\eta$ has canonical singularities and $B_\eta=0$
(see Fujino, Mori~\cite{fm, osamu}, or ~\cite{NM}).
Thus, we simply extend their methods to deal 
with the case of varieties with boundary. 

\begin{ack} I would like to thank Professors Alessio Corti,
Nick I. Shepherd-Barron and Vyacheslav V. Shokurov for 
useful discussions. This work was supported through a 
European Community Marie Curie Fellowship.
\end{ack}

%%%%%%%%%%%%%%%%%%%%%%%%%%%%%%%%%%
%%%%%%%%%%%%%%%%%%%%%%%%%%%%%%%%%%

\section{Preliminary}

%%%%%%%%%%%%%%%%%%%%%%%%%%%%%%%%%%%
%%%%%%%%%%%%%%%%%%%%%%%%%%%%%%%%%%%

We use the same notation and terminology as in
\cite{BP}. Since we use trancendental methods,
the base field is assumed to be the field of
complex numbers 
$k=\mathbb C$. However, the main results extend 
to the case when $k$ is an algebraically closed 
field of characteristic zero, by Lefschetz's 
Principle.

We collect below some results of 
Kawamata~\cite{koddim,minit}, Koll\'ar~\cite{big}
and Viehweg~\cite{wp,wpII}, with minor 
modifications.

%%%%%%%%%%%%%%%%%%%%%%%%%%%%%%%%%%
%%%%%%%%%%%%%%%%%%%%%%%%%%%%%%%%%%

\subsection{Equivariant resolutions}

%%%%%%%%%%%%%%%%%%%%%%%%%%%%%%%%%%%
%%%%%%%%%%%%%%%%%%%%%%%%%%%%%%%%%%%
Let $X$ be a complex-analytic
space which is countable at infinity, and let 
$Z\supset \Sing(X)$ be a closed complex subspace.
By Hironaka~\cite{Hir77}, there exists a proper 
morphism $\mu\colon X'\to X$ having the following 
properties:
\begin{itemize}
\item[(1)] $X'$ is smooth and $\mu$ induces an 
isomorphism $X'\setminus \mu^{-1}(Z)\to X\setminus Z$.
\item[(2)] $\mu^{-1}(Z)$ is a divisor with normal 
crossings support.
\item[(3)] If $U,V$ are open subsets of $X$ and 
$\alpha\colon U\to V$ is an isomorphism such that
$\alpha(Z|_U)=Z|_V$, then there exists a unique
isomorphism $\alpha'$ making the following diagram 
commutative:
\[ \xymatrix{
\mu^{-1}(U) \ar[r]^{\alpha'} \ar[d]_\mu & 
\mu^{-1}(V)  \ar[d]_\mu \\
 U \ar[r]^\alpha  & V  
} \]
\item[(4)] Let $\mathcal P$ be the pseudo-group of all 
the local isomorphisms $\alpha$ as in (3). Then $\mu$ 
is obtained as the composition of blowing-ups with 
closed smooth centers which are invariant under the 
natural liftings of $\mathcal P$.
\end{itemize}
We will say that $\mu$ is an {\em 
equivariant resolution of $X$ with respect to} $Z$.

\begin{lem}\label{eqvlem}
Let $X$ be a normal variety, let $R$ be a reduced
Weil divisor (possibly zero) and let $X\to S$ be a 
morphism. Let $\mu\colon X'\to X$ be an equivariant 
resolution of $X$ with respect to $\Sing(X)\cup R$, 
and let $E$ be a normal crossings divisor on $X'$ 
such that $R\subseteq \mu_*E$. Then
$$
\mu_*T_{X'/S}\langle -E\rangle=(\mu_*\Omega^1_{X'/S}
\langle E\rangle)^\vee\subseteq T_{X/S}.
$$
\end{lem}

\begin{proof} Since $X$ is normal and $\mu$ is 
birational, we have inclusions
$$
\mu_*T_{X'/S}\langle -E\rangle\subseteq
(\mu_*\Omega^1_{X'/S}\langle E\rangle)^\vee
\subseteq T_{X/S}.
$$
Fix a point $x\in X$ and let $\fa\in H^0(U,
(\mu_*\Omega^1_{X'/S}\langle E\rangle)^\vee)\subset
H^0(U,T_{X/S})$ be a vector field on an analytic 
neighborhood $U$ of $x$. There exists 
(see \cite{GK64}) a holomorphic one-parameter family 
$
\Phi\colon \{s\in {\mathbb C};|s|<\epsilon\}
\times U'\to U 
$
satisfying the following properties:
\begin{itemize}
\item[(a)] $\Phi(s_1,\Phi(s_2,x))=\Phi(s_1+s_2,x)$
for $|s_1|,|s_2|,|s_1+s_2|<\epsilon$.
\item[(b)] $\Phi(0,x)=x$ for $x\in U'$.
\item[(c)]  $\fa_x(df)=\frac{d}{ds} 
f(\Phi(s,x))|_{s=0}$ for $x\in U'$ and 
$f\in \calO_{X,x}$.
\end{itemize}
The local isomorphisms $\Phi_s$ preserve the 
singular locus of $X$. They also preserve $R$
on a big open subset of $X$, since $\mu$ is birational 
and $R\subseteq \mu_*E$. Therefore each $\Phi_s$ 
preserves $\Sing(X)\cup R$, so $\Phi$ lifts to a 
one-parameter family in an analytic neighborhood of 
$\pi^{-1}(x)$. The corresponding vector field
$\tilde{\fa}$ is tangent to the exceptional locus of
$\mu$ and to the components of $E$. Then
$\tilde{\fa}$ is a section of $T_{X'/S}\langle -E\rangle$
which lifts $\fa$, hence $\fa\in \mu_*T_{X'/S}\langle 
-E\rangle$.
\end{proof}

\begin{lem}\label{3.3}
Let $\mu\colon Y\to X$ be a birational morphism, 
let $B$ be a sheaf on $X$ and let $A$ be a torsion 
free sheaf on $Y$. Then we have an exact sequence
$$
0\to \Ext^1_X(B,\mu_*A)\to \Ext^1_Y(\mu^*B,A)\to
\Hom_X(B, R^1\mu_*A).
$$
\end{lem}

\begin{proof} If we set $\Hom_\mu(B,A)=\Hom_Y(\mu^*B,A)=
\Hom_X(B,\mu_*A)$, we have two spectral sequences 
(Ran~\cite{Ran89}):
$$
E_2^{p,q}=\Ext^p_X(B,R^q\mu_*A)
\Longrightarrow \Ext^{p+q}_\mu(B,A)
$$
$$
'E_2^{p,q}=\Ext^p_Y(L_q\mu^*B,A)
\Longrightarrow \Ext^{p+q}_\mu(B,A)
$$
The $5$-term exact sequence of edge homomorphisms
of the first spectral sequence is (denote 
$H^i=\Ext^i_\mu(B,A)$)
$$
0\to \Ext^1_X(B,\mu_*A)\to H^1\to
\Hom_X(B, R^1\mu_*A)\to \Ext^2_X(B,\mu_*A)\to
H^2.
$$
Since $\mu$ is birational, $L_1\mu^*B$ is
a torsion sheaf on $Y$.  Since $A$ is torsion free,
we obtain 
$'E_2^{0,1}=\Hom_Y(L_1\mu^*B,A)=0$.
Therefore the second spectral sequence induces an
isomorphism
$
\Ext^1_Y(\mu^*B,A)\isoto H^1.
$
\end{proof}

%%%%%%%%%%%%%%%%%%%%%%%%%%%%%%%%%%
%%%%%%%%%%%%%%%%%%%%%%%%%%%%%%%%%%

\subsection{The covering trick}

%%%%%%%%%%%%%%%%%%%%%%%%%%%%%%%%%%%
%%%%%%%%%%%%%%%%%%%%%%%%%%%%%%%%%%%

Let $f\colon X\to S$ be a smooth projective 
morphism of nonsingular varieties, and let 
$T$ be a $\Q$-divisor on $X$ such that $T\sim_\Q 0$ 
and such that the fractional part of $T$ has normal 
crossings support relative to $f$.
Let $\varphi$ be a rational function on
$X$ such that $(\varphi)+mT=0$ (and $m$
is minimal with this property). Let 
$\pi\colon \tilde{X}\to X$ be the 
normalization of $X$ in $k(X)(\sqrt[m]{\varphi})$,
and let $\nu\colon V\to \tilde{X}$ be an
equivariant resolution with respect to the singular
locus of $\tilde{X}$. We assume that the morphism
$h\colon V\to S$ is a smooth and $R^2h_*\Z_V$ has a
global section inducing a polarization on each fiber 
of $h$. 
\[ \xymatrix{
X  \ar[d]_f & \tilde{X}\ar[l]_\pi \ar[dl]_{\tilde{f}} 
& V\ar[l]_\nu\ar[dll]^h\\
S          &  &
} \]
We have an eigensheaf decomposition
$$
h_*K_{V/S}=\tilde{f}_*K_{\tilde{X}/S}=
\bigoplus_{i=0}^{m-1} 
f_*\calO_X(\lceil K_{X/S}+iT\rceil).
$$
In particular, we obtain locally free sheaves on $S$
$$
\calF^{(i)}:=f_*\calO_X(\lceil K_{X/S}+iT\rceil).
$$
The variation of polarized Hodge structure 
$(R^dh_*\C_V)_{prim}\otimes \calO_S$ 
induces semipositive Hermitian metrics on 
each $\calF^{(i)}$. Let
$$
\sigma^i_s\colon T_{S,s}\to
\Hom(H^0(X_s,\calF^{(i)}_s), 
\Ext^1_{X_s}(\Omega^1_{X_s}\langle E_s\rangle,
\calF^{(i)}_s))
$$
be the map induced by cup product with the 
Kodaira-Spencer class 
$$
\kappa_s\colon T_{S,s}\to 
\Ext^1_{X_s}(\Omega^1_{X_s}\langle E_s\rangle,
\calO_{X_s}),
$$
where $E$ is the support of the fractional 
part of $T$.

\begin{prop}\label{cc} The
following properties hold:
\begin{itemize}
\item[(1)] viewed as a subbundle of the flat
vector bundle 
$(R^dh_*\C_V)_{prim}\otimes \calO_S$, the
second fundamental forms of $\calF^{(i)}$ and
$h_*\omega_{V/S}$ are represented by $\sigma^i$
and $\oplus_{i=0}^{m-1}\sigma^i$, respectively.

\item[(2)] the Hermitian vector bundle $\calF^{(i)}$ 
is Griffiths semipositive definite, with curvature
$$
\Theta_i=-(\sigma^i)^*\wedge \sigma^i,
$$
where $(\sigma^i)^*(v)$ is the adjoint of
$\sigma^i(\bar{v})$ with respect to the
induced Hodge metrics. The curvature 
$\trace(\Theta_i)$
of
$
\det(\calF^{(i)})
$
is semipositive, and 
$\Ker(\sigma^i_s)\subset T_{S,s}$ consists
of the tangent directions along 
which $\trace(\Theta_i)$ is not positive 
definite at $s$.

\item[(3)] the maps 
$\sigma_{V_s}, \sigma_{\tilde{X}_s}, 
\oplus_{i=0}^{m-1}\sigma^i_s$ and 
$\kappa_{V_s}, \kappa_{\tilde{X}_s},
\kappa_s$ have the same kernel, respectively.

\item[(4)]
assume $S\subset \bar{S}$ is a nonsingular
compactification of $S$ such that 
$\bar{S}\setminus S$ is a simple normal crossings 
divisor and $R^dh_*\C_V$ has unipotent
local monodromies along $\bar{S}\setminus S$.
Let $\bar{\calF}^{(i)}$ be the Schmid extension of 
$\calF^{(i)}$. Then $\det(\bar{\calF}^{(i)})$ is a 
nef invertible sheaf on $\bar{S}$, and 
$$
\det(\bar{\calF}^{(i)})^{\dim (S)}
=\int_S (\frac{\sqrt{-1}}{2\pi}
\trace(\Theta_i))^{\dim(S)}.
$$
\end{itemize}
\end{prop}

\begin{proof} (1) We may shrink $S$, so that 
$V_s\to \tilde{X}_s$ is a resolution of singularities 
for every $s\in S$. We have an inclusion of Hodge 
structures $H^d(\tilde{X}_s,\C)\subset H^d(V_s,\C)$ and 
$H^{d,0}(\tilde{X}_s)= H^{d,0}(V_s)$, 
since $\tilde{X}_s$ has rational singularities 
(see~\cite{Steenbrink77}).
Therefore we have a commutative diagram
\[ \xymatrix{
T_{S,s}  \ar[d]_= \ar[r] & H^1(V_s,T_{V_s})
\ar[d]\ar[r] & 
\Hom(H^{d,0}(V_s),H^{d-1,1}(V_s)) \\
T_{S,s} \ar[r] & H^1(\tilde{X}_s,
T_{\tilde{X}_s}) \ar[r] & 
\Hom(H^{d,0}(\tilde{X}_s),H^{d-1,1}(\tilde{X}_s))
\ar[u]
} \]
where the middle row is the tangent map 
of the blow-down transformation of deformation
functors (\cite{minit}). If we identify (\cite{wpII})
$$
H^1(X_s,T_{X_s}\langle -E_s\rangle)=
H^1(\tilde{X}_s, T_{\tilde{X}_s})^G,
$$ 
we also have a commutative diagram:
\[ \xymatrix{
T_{S,s}  \ar[d]_= \ar[r] &  H^1(\tilde{X}_s,
T_{\tilde{X}_s}) \\
T_{S,s} \ar[r] & 
H^1(X_s,T_{X_s}\langle -E_s\rangle) \ar[u]
} \]
Therefore the cup product with the 
Kodaira-Spencer class preserves the
eigenspaces, so that the infinitesimal
period map
$$
\sigma_{V_s}\colon T_{S,s} \to 
\Hom(H^{d,0}(V_s),H^{d-1,1}(V_s))
$$
admits the decomposition 
$\sigma_{V_s}=\oplus_{i=0}^{m-1}\sigma^i_s$.
This is the same as the action of the Gauss-Manin 
connection on $(R^dh_*\C_V)_{prim}\otimes \calO_S$, 
hence $\sigma^i_s$ is the second fundamental 
form of $\calF^{(i)}$ in 
$(R^dh_*\C_V)_{prim}\otimes \calO_S$.

(3) Consider the infinitesimal period maps
\begin{gather*}
\sigma_{V_s} \colon T_{S,s}  
\to  \Hom(H^{d,0}(V_s),H^{d-1,1}(V_s)) \\
\sigma_{\tilde{X}_s} \colon T_{S,s} \to 
\Hom(H^{d,0}(\tilde{X}_s),H^{d-1,1}(\tilde{X}_s)) \\
\sigma^i_s\colon T_{S,s} \to \Hom(H^0(X_s,\calF^{(i)}_s), 
\Ext^1_{X_s}(\Omega^1_{X_s}\langle E_s\rangle,
\calF^{(i)}_s))
\end{gather*}
and the Kodaira-Spencer deformation classes
\begin{gather*}
\kappa_{V_s} \colon T_{S,s}\to
H^1(V_s,T_{V_s}) \\
\kappa_{\tilde{X}_s}\colon T_{S,s}\to
H^1(\tilde{X}_s,T_{\tilde{X}_s}) \\
\kappa_s \colon  T_{S,s}\to
H^1(X_s,T_{X_s}\langle -E_s\rangle)
\end{gather*}
The proof of (1) gives that $\sigma_{V_s}, 
\sigma_{\tilde{X}_s}, 
\oplus_{i=0}^{m-1}\sigma^i_s$ and 
$\kappa_{\tilde{X}_s},\kappa_s$ 
have the same kernel, respectively. 
Since $\nu$ is equivariant,
$\kappa_{V_s}$ and $\kappa_{\tilde{X}_s}$ have
the same kernel by \cite{minit}, Lemma 6.2.

(4) This follows from \cite{koddim} and \cite{big}. 
\end{proof}

%%%%%%%%%%%%%%%%%%%%%%%%%%%%%%%%%%
%%%%%%%%%%%%%%%%%%%%%%%%%%%%%%%%%%

\section{The period map}

%%%%%%%%%%%%%%%%%%%%%%%%%%%%%%%%%%%
%%%%%%%%%%%%%%%%%%%%%%%%%%%%%%%%%%%

Throughout this section, we fix the following
setup:
\begin{itemize}
\item $(X,B)$ is a log variety with Kawamata log 
terminal singularities, such that $K+B\sim_\Q 0$.
\item $f\colon X\to S$ is a projective contraction
to a nonsingular algebraic variety $S$.
\item $\mu\colon Y\to X$ is a resolution of singularities,
$K_Y+B_Y=\mu^*(K+B)$ is the log pullback, and
$E$ is the support of the fractional part of $B_Y$. 
We assume that $E$ has simple normal crossings
support and $\mu_*T_Y\langle -E\rangle$ is a reflexive
sheaf (such a resolution exists by Lemma~\ref{eqvlem}).
\item the family $(Y,E)\to S$ has relative simple
normal crossings over an open subset of $S$. Let 
$\kappa_s\colon T_{S,s}\to H^1(Y_s,T_{Y_s}\langle 
-E_s\rangle)$ be the induced Kodaira-Spencer class.
\end{itemize}

We are only interested in the induced log
variety $(X_\eta,B_\eta)$ defined over $k(S)$. 
Hence, we will shrink $S$ to a Zariski open subset 
without further notice.

\begin{prop}\label{ky}
Let $b$ be the minimal positive integer such that 
$b(K+B)\sim 0$ over the general point of $S$, and
choose a rational function $\varphi$ with 
$b(K+B)=(\varphi)$. Let $\pi\colon \tilde{Y}\to Y$ 
be the normalization of $Y$ in $k(Y)(\sqrt[b]{\varphi})$, 
and let $\nu\colon V\to \tilde{Y}$ be an equivariant 
resolution:
\[ \xymatrix{
(Y,B_Y)  \ar[d]_\mu & \tilde{Y} \ar[l] & 
 V \ar[ddll]^h \ar[l] \\
(X,B) \ar[d]_f & & \\
S          &  &
} \]
Let $\kappa_{V_s}$ be the induced Kodaira-Spencer 
class of $V\to S$, and let
$$
\sigma_{V_s}:T_{S,s}\to
\Hom(H^{d,0}(V_s), H^{d-1,1}(V_s))
$$
$$
\sigma^1_s\colon T_{S,s}\to
\Hom(H^0(Y_s,\lceil -B_{Y_s}\rceil), 
\Ext^1_{Y_s}(\Omega^1_{Y_s}\langle E_s\rangle,
\calO_{Y_s}(\lceil -B_{Y_s}\rceil)))
$$
be the maps induced by cup product with 
$\kappa_{V_s}$ and $\kappa_s$, respectively, 
where $d=\dim(X/S)$. The projective contraction
$h$ is smooth over an open subset $S^0$ of $S$.
Let $\Phi $ be the period map associated to the 
variation of polarized Hodge structure 
$(R^dh_*\C_V)_{prim}\otimes \calO_{S^0}$.

Then the maps $\kappa_s,\sigma^1_s,\kappa_{V_s},
\sigma_{V_s}$ have the same kernel for general 
$s\in S$, equal to
$
T_{\Phi ^{-1}(\Phi (s)),s}\subset T_{S,s}.
$
\end{prop}

\begin{proof} 
(1) Since $(X_s,B_{X_s})$ is a 
log variety with Kawamata log terminal singularities
for general $s$, we have an isomorphism
$$
\calO_{X_s}\isoto {\mu_s}_*
\calO_{Y_s}(\lceil -B_{Y_s}\rceil).
$$
In particular, the inclusion 
$\calO_{Y_s}\to \calO_{Y_s}(\lceil -B_{Y_s}
\rceil)$ induces an isomorphism of global sections
and identifies $\sigma^1_s$ with the composition
$$
T_{S,s}\to \Ext^1_{Y_s}
(\Omega^1_{Y_s}\langle E_s\rangle,\calO_{Y_s})
\stackrel{\lambda}{\to}
\Ext^1_{Y_s}(\Omega^1_{Y_s}\langle E_s\rangle,
\calO_{Y_s}(\lceil -B_{Y_s}\rceil)).
$$

(2) $\Ker(\lambda)\subset 
\Ker(\Ext^1_{Y_s}(\Omega^1_{Y_s}\langle E_s\rangle, 
\calO_{Y_s})\to 
\Ext^1_{X_s}({\mu_s}_*\Omega^1_{Y_s}\langle E_s\rangle, 
\calO_{X_s}))$. Indeed, 
$\calO_{Y_s}(\lceil -B_{Y_s}\rceil)$ is ${\mu_s}_*$-acyclic
by Kawamata-Viehweg vanishing. Since $X_s$ has rational
singularities, $\calO_{Y_s}$ is ${\mu_s}_*$-acyclic as well.
The natural homomorphism $\mu_s^*{\mu_s}_*
\Omega^1_{Y_s}\langle E_s\rangle\to \Omega^1_{Y_s}\langle 
E_s\rangle$ induces a commutative diagram

\[ \xymatrix{
 \Ext^1_{Y_s}(\Omega^1_{Y_s}\langle E_s\rangle, 
\calO_{Y_s}) \ar[r]\ar[d] &   
 \Ext^1_{Y_s}(\Omega^1_{Y_s}\langle E_s\rangle, 
\calO_{Y_s}(\lceil -B_{Y_s}\rceil))\ar[d] \\
 \Ext^1_{Y_s}(\mu_s^*{\mu_s}_*\Omega^1_{Y_s}\langle 
E_s\rangle, \calO_{Y_s}) \ar[r]  &   
 \Ext^1_{Y_s}(\mu_s^*{\mu_s}_*\Omega^1_{Y_s}\langle 
E_s\rangle, \calO_{Y_s}(\lceil -B_{Y_s}\rceil))
\\
 \Ext^1_{X_s}({\mu_s}_*\Omega^1_{Y_s}\langle E_s\rangle, 
\calO_{X_s}) \ar[r] \ar[u] &   
 \Ext^1_{X_s}({\mu_s}_*\Omega^1_{Y_s}\langle E_s\rangle, 
{\mu_s}_*\calO_{Y_s}(\lceil -B_{Y_s}\rceil))\ar[u]
} \]

The bottom vertical arrows are isomorphisms by 
Lemma~\ref{3.3}, and the bottom horizontal arrow is an 
isomorphism by $\calO_{X_s}={\mu_s}_*\calO_{Y_s}
(\lceil -B_{Y_s}\rceil)$. This implies the claim.

(3)
$
\Ker(\kappa_s)=\Ker(\sigma^1_s)\subseteq 
\Ker(\kappa_{X_s})
$
(the inclusion is an equality if $E$ 
is exceptional). Indeed, the inclusion 
$\Ker(\kappa_s)\subseteq\Ker(\sigma^1_s)$ is
clear. Conversely, consider a tangent vector
$
t\in \Ker(\sigma^1_s),
$
inducing a first-order infinitesimal deformation 
of $(Y_s,E_s)$:
\[ \xymatrix{
(\bar{Y}_s,\bar{E}_s)\ar[r] \ar[d]^{\bar{\mu}_s}
& (Y,E)\ar[d]^\mu \\
\bar{X}_s\ar[d]^{\bar{f}_s} \ar[r] &
X \ar[d]^f \\
\Spec k[\epsilon]\ar[r] & S
} \]
The class $\kappa_s(t)$ is represented by
the short exact sequence
$$
0\to \calO_{Y_s} \to \Omega^1_{\bar{Y}_s}\langle 
\bar{E}_s\rangle \otimes \calO_{Y_s} 
\to \Omega^1_{Y_s} \langle E_s\rangle 
\to 0.
$$
By (2), the short exact sequence
$$
0\to \calO_{X_s}\to {\mu_s}_*(\Omega^1_{\bar{Y}_s}
\langle \bar{E}_s\rangle \otimes \calO_{Y_s})\to 
{\mu_s}_*\Omega^1_{Y_s} \langle E_s\rangle \to 0
$$
admits a splitting
$
u \in \Hom_{X_s}({\mu_s}_*(\Omega^1_{\bar{Y}_s}
\langle \bar{E}_s\rangle \otimes \calO_{Y_s}), 
\calO_{X_s}).
$
It is enough to lift $u$ to $Y_s$.
We may assume 
(cf.~\cite{minit}, Lemma 6.1) that the horizontal 
arrows are injective and the vertical arrow is 
surjective in the diagram below 
\[ \xymatrix{
\calHom_{\calO_{X_s}}({\mu_s}_*(\Omega^1_{\bar{Y}_s}
\langle \bar{E}_s\rangle\otimes\calO_{Y_s}), \calO_{X_s}) 
\ar[r] &
\calHom_{\calO_{X_s}}(\mu_*\Omega^1_Y\langle E\rangle
\otimes \calO_{X_s},\calO_{X_s}) \\
\calHom_{\calO_X}(\Omega^1_X,\calO_X) &
\calHom_{\calO_X}(\mu_*\Omega^1_Y\langle E\rangle,\calO_X)
\ar[u]\ar[l]
} \]
Fix a point $x\in X_s$. Then $u$ lifts to a local 
section $\bar{u}$ of $(\mu_*\Omega^1_Y\langle E
\rangle)^\vee$ near $x$. Since $\mu_*T_Y\langle -E\rangle$
is reflexive, we have $\mu_*T_Y\langle -E\rangle=
(\mu_*\Omega^1_Y\langle E\rangle)^\vee$.
Therefore $\bar{u}$ lifts to a section of $T_Y
\langle -E\rangle$ in a neighborhood of $\mu^{-1}(x)$. 
In particular, $u$ lifts to a local splitting of 
$\kappa_s(t)$ in a neighborhood of $\mu^{-1}(x)$. 
Local liftings of $u$ are unique, hence they glue 
to a global splitting of $\kappa_s(t)$. Therefore 
$t\in \Ker(\kappa_s)$.

The inclusion $\Ker(\kappa_s)\subseteq 
\Ker(\kappa_{X_s})$ follows from the commutative diagram
\[ \xymatrix{
0\ar[r] & \calO_{X_s}\ar[r] &
{\mu_s}_*(\Omega^1_{\bar{Y}_s}
\langle \bar{E}_s\rangle \otimes \calO_{Y_s})\ar[r] &
{\mu_s}_*\Omega^1_{Y_s} \langle E_s\rangle \ar[r] & 0 \\
0\ar[r] & \calO_{X_s}\ar[r]\ar[u]^= &
\Omega^1_{\bar{X}_s} \otimes \calO_{X_s}\ar[r]\ar[u] &
\Omega^1_{X_s} \ar[r] \ar[u] & 0
} \]

(4) The natural inclusion  
$
\Ker(\kappa_s)\subseteq 
\Ker(\oplus_{i=0}^{b-1}\sigma^i_s)\subseteq\Ker(\sigma^1_s)
$
and Proposition~\ref{cc}(3) imply that the maps
$
\kappa_s, \kappa_{\tilde{Y}_s},
\kappa_{V_s}, \sigma^1_s, \sigma_{\tilde{Y}_s},
\sigma_{V_s}
$ 
have the same kernel, for general $s\in S$.
The inclusions
$$
\Ker(\kappa_{V_s})\subseteq \Ker(d\Phi _*)
\subseteq \Ker(\sigma_{V_s}),
$$
imply that the common kernel is the tangent space 
at $s$ of the fiber $\Phi ^{-1}(\Phi (s))$ of $\Phi $.
\end{proof}

\begin{thm}\label{kl}
There exist dominant morphisms $\tau\colon \bar{S}\to S$
and $\varrho\colon \bar{S}\to S^!$, with $\tau$ generically
finite and $\bar{S},S^!$ nonsingular, and there exist 
a log variety $(X^!,B^!)$ with $K_{X^!}+B^!\sim_\Q 0$ and 
a projective contraction $f^!\colon X^!\to S^!$
\[ \xymatrix{
(X,B) \ar[d]_f &  & (X^!,B^!) \ar[d]_{f^!} \\
S        & \bar{S} \ar[l]_\tau  \ar[r]^\varrho  & S^! 
} \]
satisfying the following properties:
\begin{itemize}
\item[(1)] there exists an open dense subset 
$U\subset \bar{S}$ and an isomorphism
\begin{displaymath}
\xymatrix{
(X,B) \times_S \bar{S}\vert_U \ar[rr]^\simeq \ar[dr] & &
(X^!,B^!)\times_{S^!} \bar{S}\vert_U \ar[dl]\\
   &  U  &                                    }
\end{displaymath}
\item[(2)] $\kappa_{s^!}$ is injective for general
points $s^!\in S^!$.
\end{itemize}
\end{thm}

\begin{proof} We may assume that $h\colon V\to S$ 
is a smooth projective morphism. Let $S\subset \bar{S}$ 
be a nonsingular compactification with simple normal 
crossing boundary.
By~\cite{Griffiths}, the period map $\Phi $ of the variation 
of polarized Hodge structure 
$(R^dh_*{\mathbb C}_V)_{prim}\otimes \calO_S$ extends to
a proper analytic morphism defined on an open subvariety 
of $\bar{S}$. By~\cite{koddim}, $\Phi$ is 
bimeromorphic with a rational map
$$
\Phi \colon S \to S^!
$$
and we may further shrink $S$ so that $\Phi$ is regular 
everywhere on $S$. After we replace $S$ by an \'etale 
open set, the following properties hold:
\begin{itemize}
\item[(a)] $(Y,E)\to S$ is a smooth relative pair and
${\mu_s}_*\calO_{Y_s}=\calO_{X_s}$ for every $s\in S$.
\item[(b)] $\Ext^1_{Y_s}(\Omega^1_{Y_s}\langle E_s\rangle,
\calO_{Y_s})$ and $\Ext^1_{X_s}(\Omega^1_{X_s},\calO_{X_s})$
have constant dimension for $s\in S$.
\item[(c)] $\Ker(\kappa_s)=T_{\Phi^{-1}(\Phi(s))}$ for every 
$s\in S$.
\item[(d)] $\Phi\colon S\to S^!$ has a section
$i\colon S^!\to S$.
\end{itemize}
By base change via the section $i$, we induce a 
family on $S^!$: 
\[ \xymatrix{
(Y,B_Y) \ar[d]^\mu & (Y^!,B_{Y^!}) \ar[d]^{\mu^!} \\
(X,B) \ar[d]^f & (X^!,B^!)\ar[d]^{f^!} \\
S\ar[r]^\Phi & S^!
} \]
According to Proposition~\ref{ky} and its proof, 
the families $(Y,E)\to S$ and $X\to S$ are first-order 
infinitesimally trivial when restricted to the fibers of 
$\Phi$. The same holds for the map $Y\to X$ over $S$, 
by (a) (see~\cite{Ran89}). By ~\cite{Gr74} and 
~\cite{Ka78}, the family
$$
(Y,E)\to X\to S
$$
is locally trivial when restricted to the fibers of 
$\Phi$.
Consider the subfunctor of
$$
\Isom_S(X,X^!\times_{S^!} S)\times 
\Isom_S((Y,E),(Y^!,E^!)\times_{S^!} S),
$$
making the obvious diagrams commutative. This subfunctor
is representable by a scheme $I/S$, and the map $I\to S$ 
is surjective from the above considerations. After replacing 
$S$ (and $S^!$ accordingly) by an \'etale open subset, 
we may assume that $I/S$ has a section.
Consequently, we have global isomorphisms
\[ \xymatrix{
(Y,E)\ar[r] \ar[d] & (Y^!,E^!)\times_{S^!} S \ar[d] \\
X\ar[d] \ar[r]     & X^!\times_{S^!} S \ar[d] \\
S\ar[r]^= & S
} \]
Since $(X,B)$ is klt and $B$ is effective, we 
have 
$B=\mu_*E$ and 
$B^!\times_{S^!} S=(\mu^!\times_{S^!} 1_S)_*
(E^!\times_{S^!} S)$.
Therefore $X\to X^!\times_{S^!} S$ is in fact 
an isomorphism of pairs over $S$
$$
(X,B)\isoto (X^!,B^!)\times_{S^!} S.
$$
By (c), the Kodaira-Spencer class 
$\kappa_{s^!}$ is injective for $s^!\in S^!$. 
\end{proof}

%%%%%%%%%%%%%%%%%%%%%%%%%%%%%%%%%%
%%%%%%%%%%%%%%%%%%%%%%%%%%%%%%%%%%

\section{Lc-trivial fibrations}

%%%%%%%%%%%%%%%%%%%%%%%%%%%%%%%%%%%
%%%%%%%%%%%%%%%%%%%%%%%%%%%%%%%%%%%

Recall (\cite{BP}) that an {\em lc-trivial 
fibration} $f\colon (X,B)\to Y$ consists of 
a contraction $f\colon X\to Y$ of proper normal 
varieties and a log pair $(X,B)$, subject to the 
following conditions:
\begin{itemize}
\item[(1)] $(X,B)$ has Kawamata log terminal 
singularities over the generic point of $Y$,
\item[(2)] $\rank f_*\calO_X(\lceil \bA(X,B)
\rceil)=1$,
\item[(3)] there exist a positive integer $r$, a
rational function $\varphi\in k(X)^\times$ and a 
$\Q$-Cartier divisor $D$ on $Y$ such that
$$
K+B+\frac{1}{r}(\varphi)=f^*D.
$$
\end{itemize}

The lc-trivial fibration $f\colon (X,B)\to Y$
induced $\Q$-b-divisors $\bB$ and $\bM$ of $Y$,
called the {\em discriminant} and {\em moduli} 
$\Q$-b-divisor respectively. By \cite{BP},
Theorem 0.2, the moduli $\Q$-b-divisor is
b-nef. This means that there exists a proper
birational model $Y'$ of $Y$ such that 
$\bM=\overline{\bM_{Y'}}$ and $\bM_{Y'}$ 
is a nef $\Q$-Cartier divisor. 

\begin{prop}\label{bc} 
Let $f\colon (X,B)\to Y$ be an
lc-trivial fibration. Let $\varrho\colon Y'\to Y$
be a surjective morphism from a proper
normal variety $Y'$, and let 
$f'\colon (X',B_{X'})\to Y'$ be an lc-trivial
fibration induced by base change
\[ \xymatrix{
(X,B) \ar[d]_f & (X',B_{X'}) \ar[l]_{\varrho_X}
\ar[d]_{f'} \\
Y        & Y'\ar[l]_\varrho
} \]
Let $\bM$ and $\bM'$ be the corresponding moduli 
$\Q$-b-divisors. Then $$\varrho^*\bM=\bM'.$$
\end{prop}

\begin{proof} (1) To define the induced lc-trivial 
fibration, we may assume that $X'$ is the 
normalization of the main component of $X\times_Y Y'$. 
Also, we may assume that $X$ is nonsingular. Then
there exists a canonical divisor $K_{X'}$ on $X'$
such that $A=K_{X'}-\varrho_X^*K$ is an 
$f'$-vertical Weil divisor. 
Define $B_{X'}=\varrho_X^*B-A$, so that
$$
K_{X'}+B_{X'}+\frac{1}{r}(\varrho_X^*\varphi)=
{f'}^*(\varrho^*D).
$$
The moduli $\Q$-b-divisor $\bM'$ of the lc-trivial
fibration $(X', B_{X'})$ is independent of the 
above choice of $K_{X'}$, by~\cite{BP}, Remark 3.3. 

(2) $\varrho^*\bM=\bM'$ if $\varrho$ 
is generically finite. This is clear if $\varrho$ 
is birational. If 
$\varrho$ is a finite morphism, $\varrho^*\bM=\bM'$ 
by \cite{thesis}, Theorem 3.2 and \cite{BP}, 
Remark 3.3. 

(3) The general case follows from (1) and the 
compatibility with base change of the canonical 
extension of a variation of Hodge structure with 
unipotent local monodromies at infinity (same 
argument as in the proof of \cite{BP}, Theorem 2.7).
\end{proof}

\begin{defn} A $\Q$-b-divisor $\bM$ of $Y$ is 
called {\em b-nef and good} if there exists a 
proper birational model $Y'$ of $Y$, endowed with a 
proper contraction $h\colon Y'\to Z$, such that
\begin{itemize}
\item[(1)] $\bM_{Y'}\sim_\Q h^*H$, for some 
nef and big $\Q$-divisor $H$ of $Z$.
\item[(2)] $\bM=\overline{\bM_{Y'}}$.
\end{itemize}
\end{defn}

\begin{thm}\label{klglobal} 
Let $f\colon (X,B)\to Y$ be an lc-trivial fibration 
such that the geometric generic fiber 
$X_{\overline{\eta}}=X\times_Y \Spec(\overline{k(Y)})$ 
is a projective variety and $B_{\overline{\eta}}$ is 
effective.
Then there exists a diagram
\[ \xymatrix{
(X,B) \ar[d]_f &  & (X^!,B^!) \ar[d]_{f^!} \\
Y        & \bar{Y} \ar[l]_\tau  \ar[r]^\varrho  & Y^! 
} \]
satisfying the following properties:
\begin{itemize}
\item $f^!\colon (X^!,B^!)\to Y^!$ is an 
lc-trivial fibration.
\item $\tau$ is generically finite and 
surjective and $\varrho$ is surjective.
\item there exists a nonempty open subset 
$U\subset \bar{Y}$ and an isomorphism
\[ \xymatrix{
(X,B) \times_Y \bar{Y}|_U \ar[dr] \ar[rr]^\simeq
&  & (X^!,B^!)\times_{Y^!} \bar{Y}|_U \ar[dl] \\
        & U  &  
} \]
\item Let $\bM, \bM^!$ be the corresponding moduli 
$\Q$-b-divisors. Then $\bM^!$ is b-nef and 
big and $\tau^*\bM=\varrho^*(\bM^!)$.
\end{itemize}
In particular, the moduli $\Q$-b-divisor $\bM$
is b-nef and good.
\end{thm}

\begin{proof} Note that the restriction of
$f\colon (X,B)\to Y$ to an appropriate open 
subset of $Y$ satisfies the assumptions of
Section 2.

(1) Assume that $\kappa_s$ is injective for
sufficiently general points $s\in Y$. Then 
$\bM$ is b-nef and big.

Indeed, let $\mu\colon X'\to X$ be an equivariant 
resolution of $X$ with respect to 
$\Sing(X)\cup \Supp(B)$. Let $\tilde{X'}\to X'$ be 
the normalization of $X'$ in $k(X')(\sqrt[b]{\varphi})$ 
and let $V\to \tilde{X'}$ be an equivariant resolution
with respect to the singular locus of $\tilde{X'}$, and
let $f'\colon X'\to Y$ be the induced morphism:
\[ \xymatrix{
X'  \ar[d]_\mu & \tilde{X'}\ar[l] & 
V\ar[l]  \ar[ddll]^h\\
X\ar[d]_f &   & \\
Y          &  &
} \]
After a generically finite base change of $Y$,
we may assume that $\bM$ descends to $Y$ and 
$V/Y$ is birational to $V'/Y$ which has simple 
normal crossings degeneration and is semi-stable 
in codimension one. Under these assumptions, we
infer by \cite{BP}, Lemma 5.2 that $\bM_Y$ is a
nef Cartier divisor and $\calO_Y(\bM_Y)$ is 
isomorphic to the Schmid extension
$$
\bar{\calF}^{(1)}=
f'_*\calO_{X'}(\lceil -B_{X'}+{f'}^*B_Y+{f'}^*\bM_Y\rceil).
$$
By Proposition~\ref{ky}, $\sigma^1_s$ and 
$\kappa_s$ have the same kernel, hence $\sigma^1_s$ 
is injective. By Proposition~\ref{cc}, the 
invertible sheaf $\bar{\calF}^{(1)}$ has positive
self-intersection. Therefore $\bM_Y$ is a nef and 
big Cartier divisor, hence $\bM=\overline{\bM_Y}$ 
is b-nef and big.

(2) We may assume that the base space $Y$ is 
nonsingular. By Theorem~\ref{kl}, there exists a 
diagram
\[ \xymatrix{
(X,B) \ar[d]_f &  & (X^!,B^!) \ar[d]_{f^!} \\
Y        & \bar{Y} \ar[l]_\tau  \ar[r]^\varrho  & Y^! 
} \]
satisfying the following properties:
\begin{itemize}
\item[(a)] $\tau$ is generically finite and 
surjective and $\varrho$ is surjective.
\item[(b)] $f^!\colon (X^!,B^!)\to Y^!$ is an 
lc-trivial fibration.
\item[(c)] there exists a nonempty open subset 
$U\subset \bar{Y}$ and an isomorphism
\[ \xymatrix{
(X,B) \times_Y \bar{Y}|_U \ar[dr] \ar[rr]^\simeq
&  & (X^!,B^!)\times_{Y^!} \bar{Y}|_U \ar[dl] \\
        & U  &  
} \]
\end{itemize}
The moduli $\Q$-b-divisor $\bM^!$ is b-nef
and big, by (1). 
Let $\bar{f}\colon \bar{X}\to \bar{Y}$ be a 
contraction which is birationally induced by 
both $f$ and $f^!$, and let $(\bar{X}/\bar{Y},
B_{\bar{X}})$ and $(\bar{X}/\bar{Y},B^!_{\bar{X}})$ 
be the lc-trivial fibrations induced by base change. 
From (c), there exists a $\Q$-Cartier divisor $L$ on 
$\bar{Y}$ such that $B^!_{\bar{X}}=B_{\bar{X}}+
\bar{f}^*L$. Therefore the lc-trivial fibrations 
$(\bar{X}/\bar{Y},B_{\bar{X}})$ and $(\bar{X}/\bar{Y},
B^!_{\bar{X}})$ have the same moduli $\Q$-b-divisor 
(\cite{BP}, Remark 3.3). We conclude by 
Proposition~\ref{bc} that
$
\tau^*(\bM)=\varrho^*(\bM^!).
$
\end{proof}

\begin{prop}\label{num} 
Let $f\colon X\to Y$ be a morphism of projective 
manifolds, and let $\Sigma_Y\subset Y$ be a simple 
normal crossings divisor such that
\begin{itemize}
\item[(i)] $f$ is smooth over $Y\setminus \Sigma_Y$.
\item[(ii)] $f$ is semi-stable over codimension one 
points of $Y$.
\end{itemize}
Let $\calL\subset f_*\omega_{X/Y}$ be a direct summand
invertible sheaf such that $\calL\equiv 0$. Then 
$\calL^{\otimes r}\simeq \calO_X$ for some positive 
integer $r$.
\end{prop}

\begin{proof} 
Consider the variation of polarized Hodge structure
$$
H=R^df_*{\mathbb C}_{X^0}\otimes \calO_{Y^0},
$$ 
where
$Y^0=Y\setminus \Sigma_Y$ and $d=\dim(X/Y)$. 
The sheaf 
$$
\calL|_{Y^0} \subset f_*\omega_{X/Y}|_{Y^0}=F^dH
$$
has an induced Hermitian metric $h$ with semi-positive
curvature form $\Theta$. We claim that the curvature
is trivial. Indeed, let $v\in T_{Y,s}$ be a tangent
vector at a point $s\in Y^0$. There exists a projective
curve $C\subset Y$ such that $T_{C,s}={\mathbb C}v$, and
let $\nu\colon C^\nu \to C$ be its normalization. 
Since the local monodromies of $H$ at infinity are 
unipotent by assumption, we infer by Kawamata~\cite{Ka81}
that
$$
\deg(\calL|_C)=\frac{\sqrt{-1}}{2\pi}\int_{\nu^{-1}(Y^0)}
\nu^*\Theta.
$$
Since $\calL$ is numerically trivial, we have 
$\nu^*\Theta=0$. Since $\nu$ is an isomorphism near
$s$, we conclude that $v$ lies is a null direction 
of $\Theta$.

Therefore $\Theta$ is trivial, i.e., 
$\calL|_{Y^0}\subset H$ is a local 
subsystem. By Deligne, there exists a positive
integer $r$ such that $(\calL|_{Y^0})^{\otimes r}$
is a trivial local system. By Kawamata, 
$f_*\omega_{X/Y}$ is the canonical extension of
$f_*\omega_{X/Y}|_{Y^0}$. The same property holds 
for its direct summand $\calL$. Since the local
monodromies are unipotent, the canonical extension
commutes with tensor products. Therefore
$\calL^{\otimes r}\simeq \calO_X$.
\end{proof}

\begin{thm}\label{torsion} 
Let $f\colon (X,B)\to Y$ be an lc-trivial fibration. 
If $\bM$ is b-numerically trivial, then $\bM\sim_\Q 0$.
\end{thm}

\begin{proof} This is similar to the proof of 
\cite{BP}, Theorem 0.1.
After a finite base change~\cite{BP}, Lemma 5.1, we 
may assume that the induced root fiber space 
$h\colon V\to Y$ is semi-stable in codimension one
and $\bM$ descends to $Y$. 
By construction, the invertible sheaf
$\calO_Y(\bM_Y)\subset h_*\omega_{V/Y}$ is a 
direct summand, and $\bM_Y\equiv 0$. 
 
We conclude by Proposition~\ref{num} that $\bM_Y$ is 
torsion. Therefore $\bM$ is torsion.
\end{proof}

%%%%%%%%%%%%%%%%%%%%%%%%%%%%%%%%%%
%%%%%%%%%%%%%%%%%%%%%%%%%%%%%%%%%%

\section{Applications}

%%%%%%%%%%%%%%%%%%%%%%%%%%%%%%%%%%%
%%%%%%%%%%%%%%%%%%%%%%%%%%%%%%%%%%%

\begin{thm} Let $(X,B)$ be a projective log variety 
with Kawamata log terminal singularities, let
$f\colon X\to Y$ be a contraction to a proper
normal variety $Y$ and let $\omega$ be a $\Q$-Cartier 
divisor on $Y$ such that
$$
K+B\sim_\Q f^*\omega.
$$
Then there exists a $\Q$-Weil divisor $B_Y$ such 
that $(Y,B_Y)$ is a log variety with Kawamata log
terminal singularities and 
$
\omega \sim_\Q K_Y+B_Y.
$
\end{thm}

\begin{proof} We may write 
$K+B+\frac{1}{r}(\varphi)=f^*\omega$,
where $r$ is a positive integer and $\varphi$
is a rational function on $X$. Thus, 
$f\colon (X,B)\to Y$ is an lc-trivial 
fibration. Denote by $\bB$ and $\bM$ the
induced different and moduli $\Q$-b-divisors.
The assumptions of Theorem~\ref{klglobal} are
satisfied, so we may find
a high resolution $\sigma\colon Y'\to Y$ 
such that $\bM=\overline{\bM_{Y'}}$ and
$$\bM_{Y'}\sim_\Q h^*A,$$
where $h\colon Y'\to Z$
is a contraction to a normal projective variety 
$Z$ and $A$ is a nef and big divisor on $Z$.
Consider the lc-trivial fibration induced by
base change with $\sigma$:
\[ \xymatrix{
(X,B) \ar[d]_f & (X',B_{X'})\ar[l]_\nu
\ar[d]_{f'}\\
Y         & Y' \ar[l]_{\sigma}
} \]
Then $\sigma^*\omega=K_{Y'}+\bB_{Y'}+\bM_{Y'}$.
Since $\bM$ descends to $Y'$, Inversion of
Adjunction holds for $f'$ (\cite{BP}, Theorem 3.1).
Therefore $(Y',\bB_{Y'})$ is a log pair with 
Kawamata log terminal singularities. We may
find an effective $\Q$-divisor $E$ on $Y'$
such that $\bM_{Y'}\sim_\Q E$ and $(Y',\bB_{Y'}+E)$
has Kawamata log terminal singularities. If we
set $B_Y=\sigma_*(\bB_{Y'}+E)$, then
$$
\sigma^*(K_Y+B_Y)=K_{Y'}+\bB_{Y'}+E\sim \sigma^*\omega,
$$
and $(Y,B_Y)$ is a log pair with Kawamata log terminal 
singularities.
Since $B$ is effective, $\bB_Y=\sigma_*\bB_{Y'}$
is effective. Therefore $B_Y$ is effective, i.e.,
$(Y,B_Y)$ is a log variety.
\end{proof}

\begin{thm}\label{AB0} 
Let $(X,B)$ be a projective log variety 
with Kawamata log terminal singularities such that
$K+B\equiv 0$. Then $K+B\sim_\Q 0$.
\end{thm}

\begin{proof}
The variety $X$ has rational singularities
since $(X,B)$ has Kawamata log terminal 
singularities. Therefore the Albanese map of 
$X$ is a morphism (\cite{minit}). 

We use induction on $\dim(X)$.
If $q(X)=0$, the numerically trivial divisor 
$K+B$ is certainly a torsion divisor. Assume now 
that $q(X)>0$.
Let $f\colon X\to Y$ be the Stein factorization 
of the Albanese map of $X$. 
The geometric generic fibre $(X_{\bar{\eta}},
B_{\bar{\eta}})$ satisfies the
same properties as $(X,B)$ and $\dim X_{\bar{\eta}}<
\dim X$. Therefore $K_{X_{\bar{\eta}}}+B_{\bar{\eta}}
\sim_\Q 0$, by induction. In particular, $K+B$ is 
numerically trivial and $\Q$-linearly 
equivalent to an $f$-vertical divisor.

Therefore we may choose a sufficiently
high resolution $\mu\colon Y'\to Y$ and a 
diagram induced by base change
\[ \xymatrix{
(X,B) \ar[d]_f & (X',B_{X'})\ar[l]_\nu
\ar[d]_{f'}\\
Y          & Y' \ar[l]
} \]
such that $f'\colon (X',B_{X'})\to Y'$
is an lc-trivial fibration and every prime
divisor on $X'$ is exceptional on $X$ if it 
is exceptional on $Y'$. We have
$$
\nu^*(K+B)+{f'}^*(\bB_{Y'}^-)=
{f'}^*(K_{Y'}+\bB_{Y'}^++\bM_{Y'}),
$$
where $\bB_{Y'}=\bB_{Y'}^+- \bB_{Y'}^-$ is 
the decomposition of the discriminant on $Y'$ 
into positive and negative components. 
The effective $\Q$-divisor ${f'}^*(\bB_{Y'}^-)$ 
is exceptional on $X$, since it is supported by
the negative part of $B_{X'}$ and by 
$f'$-exceptional divisors.

By \cite{Ueno}, Theorem 10.3, $\kappa(Y',K_{Y'})\ge 0$. 
Since $\bM$ is b-nef and good, we also have 
$\kappa(Y',\bM_{Y'})\ge 0$. 
Therefore $\kappa(X,K+B)\ge 0$, hence
$K+B\sim_\Q 0$.
\end{proof}

\begin{thm}\label{ab} 
Let $(X,B)$ be a projective log variety with 
Kawamata log terminal singularities, such that the 
log canonical class $K+B$ is nef, of nef dimension
$n$. Assume that the Log Minimal Model Program and
the Log Abundance Conjecture are valid in dimension
$n$. 

Then the linear system $|k(K+B)|$ is base point 
free for large and divisible integers $k$. 
\end{thm}

\begin{proof} The argument of \cite{NM}, Theorem 5.1
is still valid, provided we replace
the references \cite{minit}, Theorem 8.2 and
\cite{NM}, Theorem 4.5 with Theorem~\ref{AB0} and
Theorem~\ref{klglobal}, respectively.
\end{proof}

For the rest of this section, we generalize some
results of Viehweg~\cite{wp} and Kawamata~\cite{minit} 
on Ueno's Conjecture K~\cite{Ueno}.

\begin{prop}\label{extend} 
Let $f\colon (X,B)\to Y$ be an lc-trivial 
fibration of normal projective varieties such 
that there exists an isomorphism 
$$
\Phi\colon (X, B)|_U\isoto 
(F,B_F)\times U\mbox{ over } U,
$$ 
where $U\subset X$ is a non-empty open subset.
Let 
$$
Y^0=Y\setminus \{\Sing(Y)\cup \Supp(\bB_Y)\cup
f(\Supp(B_-^v)\},
$$
where $B_-^v$ is the negative and vertical part
of $B$. Then $\Phi$ extends to an isomorphism
$\Phi\colon (X,B)|_{Y^0} \isoto 
(F, B_F)\times Y^0$ over $Y^0$.
\end{prop}

\begin{proof}
(1) Let $E$ be a prime divisor on $X$ such 
that $\codim_X(f(E))\ge 2$. Then 
$f(E)\cap Y^0=\emptyset$.

Indeed, $\bM\sim_\Q 0$ by (ii). In particular, 
$\bM$ descends to $Y$, i.e., Inversion of Adjunction 
(\cite{BP}, Theorem 3.1) holds for $f$. 
If $B$ is negative at $E$, then 
$f(E)\cap Y^0=\emptyset$. Otherwise, the log 
discrepancy of $(X,B)$ at the generic point of $E$ is 
$1-\mult_E(B)\le 1$.
By Inversion of Adjunction, the minimal log discrepancy 
of $(Y,\bB_Y)$ at the generic point of $f(E)$ is at 
most one. 
Since $B$ is effective, $\bB_Y$ is effective. Since
$f(E)$ is a subvariety of codimension at least two, 
it must be contained either in the singular locus
of $Y$, or in the support of $\bB_Y$.

(2) The rational map $\Phi\colon X|_{Y^0}\dashrightarrow
F\times Y^0$, and its inverse, are isomorphisms in
codimension one. Moreover, $B|_{Y^0}$ is a horizontal
divisor.

By assumption, $\Phi$ is an isomorphism above 
the generic point of $Y^0$. By (1), a prime divisor 
on $X|_{Y^0}$ is vertical over $Y^0$ if and only if it
maps onto a prime divisor of $Y^0$. Therefore we may 
assume that $Y$ is a curve and it suffices to show that
$\Phi$ induces a birational map $X_P\dashrightarrow
F\times P$ for every point $P\in Y\setminus \Supp(\bB_Y)$.

We have 
$
\bA(X,B+f^*P)=\bA(F\times Y, B_F\times Y+F\times 
\bB_Y+F\times P).
$
The log variety
$$
(F\times Y, B_F\times Y+F\times P)=
(F,B_F)\times (Y,P)
$$
has log canonical singularities near $F\times P$, 
with $F\times P$ the unique lc place over $P$.
Since $\bB_Y=0$ near $P$, $(X,B+f^*P)$ has log canonical 
singularities near $X_P$, with $F\times P$ the unique
lc place over $P$. Since $B$ is effective at the 
components of $X_P$, this implies that $X_P$ is a 
reduced prime divisor and $\Phi$ induces a birational 
map $X_P\dashrightarrow F\times P$. Moreover,
$B$ has multiplicity zero at $X_P$.

(3) We have an induced isomorphism 
$\Phi\colon V_1\isoto V_2$, where $V_1, V_2$ are
big open subsets of $X|_{Y^0}$ and $F\times Y^0$, 
respectively. Fix a point on $Y^0$ with a local
chart $(\Delta;t_1, \ldots, t_n)$. The sections
$p^*\frac{\partial}{\partial t_i}\in 
H^0(F\times\Delta,(p^*\Omega^1_Y)^\vee)$ 
lift to vector fields
$$
\fa_i\in H^0(F\times \Delta, T_{F\times Y}).
$$
The sheaves $T_X$, $T_{F\times Y}$, $(p^*\Omega^1_Y)
^\vee$ and $(f^*\Omega^1_Y)^\vee$ are reflexive,
hence $\Phi$ induces vector fields
$$
\fa'_i\in H^0(X|_\Delta, T_{X})
$$
which lift $f^*\frac{\partial}{\partial t_i}\in 
H^0(X|_\Delta, (f^*\Omega^1_Y)^\vee)$. 
By \cite{GK64}, these vector fields 
define one-parameter groups of automorphisms 
in a neighborhood of the fixed fibre, which 
in turn define a trivialization of $f$ near 
the fixed fibre.
By construction, the one-parameter groups
of $\fa_i$ and $\fa'_i$ are compatible via
$\Phi$. Therefore $\Phi$ is an isomorphism
near the fixed fibre.

As for the boundary, note that $B|_{Y^0}$ is
horizontal over $Y^0$. Since $\Phi$ preserves
the boundaries over the generic point of $Y^0$,
the isomorphism $\Phi\colon X|_{Y^0}\isoto 
F\times Y^0$ satisfies $\Phi(B)=B_F\times Y^0$.
\end{proof}

\begin{exmp} Let $Y$ be a surface with DuVal 
singularities.
Then the minimal resolution $f\colon X\to Y$ is an lc 
trivial fibration with $\bB_Y=0$ and $\bM=0$. It is an 
isomorphism only outside the DuVal singularities.
\end{exmp}

\begin{prop}\label{aut} 
Let $(X,B)$ be a proper log variety 
with Kawamata log terminal singularities, such that 
$\kappa(X,K+B)\ge 0$. Let $\Aut^0(X,B)$
be the connected component of the group scheme
$$
\Aut(X,B)=\{\sigma\in \Aut(X); \sigma_*(B)=B\}
$$
Then $\Aut^0(X,B)$ is an Abelian variety.
\end{prop}

\begin{proof} It is known that $\Aut^0(X)$ is
an algebraic group~\cite{FGA, MO}, so
its closed subgroup $\Aut^0(X,B)$ is an algebraic 
group. Assume by contradiction that $\Aut^0(X,B)$ 
contains a linear algebraic group. Then $\Aut^0(X,B)$ 
contains a connected one dimensional linear group 
$G=\G_m$ or $\G_a$ (\cite{Ros56}).

The closed subset $\Sing(X)\cup \Supp(B)$ is
$G$-invariant. By \cite{Ros56}, Theorem 10, 
there exists a $G$-invariant open subset 
$U\subset X\setminus (\Sing(X)\cup \Supp(B))$
and a $G$-invariant isomorphism
$$
U\isoto G\times V,
$$
where $G$ acts on $G\times V$ only on the first 
factor, by translations.
In particular, $V$ is non-singular. Choose a 
compactification $V\subset Y$ such that 
$Y\setminus V$ is a simple normal crossings 
divisor. By Hironaka's resolution of singularities,
there exists a diagram
\[ \xymatrix{
    & X'\ar[dl]_f\ar[dr]^g & \\
X   &                      &  {\mathbb P}^1\times Y
} \]
such that $f$ is an isomorphism over $U$, $g$ is
an isomorphism over $G\times V$, $X'$ is proper and
nonsingular and 
$$
f^{-1}(X\setminus U)=g^{-1}({\mathbb P}^1\times Y
\setminus G\times V)=\sum E_i
$$
is a simple normal crossings divisor on $X'$. Let
$f^*(K+B)=K+B_Y$, and let $B^+_Y$ be the positive
part of $B_Y$. We have $\lfloor B^+_Y\rfloor=0$, 
since $(X,B)$ has Kawamata log terminal singularities.
Since $B$ is effective, we have
$$
\kappa(X', K_{X'}+B_Y^+)=\kappa(X,K+B)\ge 0.
$$
In particular, $\kappa({\mathbb P}^1\times Y,
g_*(B_Y^+))\ge 0$. But
$$
({\mathbb P}^1\times Y, g_*(B_Y^+))=
({\mathbb P}^1,B_{{\mathbb P}^1})\times
(Y, B_Y),
$$
where $B_{{\mathbb P}^1}$ and $B_Y$ are boundaries
supported by ${\mathbb P}^1\setminus G$ and
$Y\setminus V$, respectively. In particular,
$\kappa({\mathbb P}^1,K_{\mathbb P^1}+
B_{{\mathbb P}^1})\ge 0$. Since $B_{{\mathbb P}^1}$
is a boundary, this implies that $G=\G_m$ and
$B_{\mathbb P^1}$ is the reduced sum of two points. 
This contradicts $\lfloor B^+_Y\rfloor=0$.
\end{proof}

\begin{thm}\label{str}
Let $f\colon (X,B)\to Y$ be an lc-trivial fibration 
such that $X_{\bar{\eta}}$ is projective,
$B$ is effective over a big open subset of $Y$, 
$\bB_Y=0$ and $\bM\sim_\Q 0$.
Then there exists a finite Galois covering 
$\tau\colon Y'\to Y$ such that:
\begin{itemize}
\item[(i)] $\tau$ is \'etale in codimension one.
\item[(iii)] there exists a non-empty open subset 
$U\subset Y'$ and an isomorphism $(X,B)\times_ Y Y'|_U
\isoto (F,B_F)\times Y'|_U$ over $U$. 
\end{itemize}
\end{thm}

\begin{proof} (1) The fibers of $f$ are reduced over 
a big open subset of $Y$.
Indeed, we may assume that $Y$ is a curve. Let $P\in Y$
be a point and let $f^*P=\sum_i m_i E_i$. We have 
$$
1=1-\mult_P(\bB_Y)\le \min_i\frac{1-\mult_{E_i}(B)}{m_i}
\le \frac{1}{m_i}.
$$
Therefore the fiber $f^*P$ is reduced. 
In particular, there exists a big open subset 
$Y^0\subset Y$ such that $B$ is horizontal over 
$Y^0$ and $f\colon X|_{Y^0}\to Y^0$ is smooth on 
a big open subset of $X|_{Y^0}$.

(2) Since $\bM\sim_\Q 0$, there exist by 
Theorem~\ref{klglobal} a generically
finite morphism $\tau\colon W\to Y$ from a nonsingular 
proper variety $W$ and a nonempty open subset $U\subseteq W$ 
such that $(X,B)\times_Y W|_U$ and $(F,B_F)\times W|_U$ 
are isomorphic over $U$. We may assume that the field extension
$k(W)/k(Y)$ is Galois and $G=\Gal(k(W)/k(Y))$ acts
regularly on $W$. After possibly shrinking $Y^0$, we
may assume that $W^0:=\tau^{-1}(Y^0)\to Y^0$ is 
a finite, flat Galois covering.
Let $X'\to X\times_Y W$ be the normalization morphism: 
\[ \xymatrix{
X \ar[d]_f & X' \ar[l]_{\tau'}\ar[d]^{g} \\
Y          &  W\ar[l]_\tau
} \]

We claim that $X'\to X\times_Y W$ is an isomorphism
above $W^0$.
Indeed, restricted to $Y^0$, $f$ is smooth on a big open 
subset and $\tau'$ is finite, hence $X\times_Y W|_{W^0}\to
W^0$ is smooth on a big open subset of $X\times_Y W|_{W^0}$. 
Since $W$ is nonsingular, the singular locus 
of $X\times_Y W|_{W^0}$ has codimension at least two. 
Furthermore, $X\times_Y W|_{W^0}$ is $S_2$ since 
$X\times_Y W|_{W^0}\to X|_{Y^0}$ is finite and flat, and
$X$ is $S_2$ (\cite{Matsumura}, 21.B Theorem 50). 

(3) By finite base change, $f\colon (X,B)\to Y$ induces
an lc-trivial fibration $g\colon (X', B_{X'}) \to W$, with 
discriminant $\bB_W$ and moduli b-divisor $\bM'\sim_\Q 
\tau^*\bM \sim_\Q 0$. Set $B':=B_{X'}-g^*(\bB_W)$.
From above, we obtain 
$$
g(\Supp(B'))\cap W^0=\emptyset\mbox{ and } \bB'_W=0.
$$
The lc-trivial fibration $g\colon (X',B')\to W$
satisfies the assumptions of Proposition~\ref{extend},
which gives an isomorphism over $W^0$
\[ \xymatrix{
(X',B_{X'})|_{W^0} \ar[dr]_g \ar[rr]^\Phi &  & (F,B_F)
\times W^0 \ar[dl]^p\\
        &  W^0 &
} \]
The Galois group $G=\Gal(W^0/Y^0)$ acts regularly
on $W^0$ and $X'|_{W^0}=X|_{Y^0}\times_{Y^0} W^0$, 
and $g$ is $G$-invariant. We have an induced action
of $G$ on $(F,B_F)\times W^0$ so that $\Phi$ is a 
$G$-invariant isomorphism.

(4) Let $H$ be the subgroup of $G$ generated
by the ramification groups $I(P)$, for every
prime divisor $P$ of $W^0$. Then 
$$
\sigma'=\mbox{id}_F\times \sigma \mbox{ for }
\sigma \in H,
$$
where $\sigma'$ is the automorphism of
$(F,B_F)\times W^0$ induced by $\sigma$.
Indeed,
$\sigma'\circ (\mbox{id}_F\times \sigma)^{-1}$ is an
automorphism of $(F,B_F)\times W^0$ over $W^0$, 
inducing a morphism
$
s_\sigma\colon W^0\to \Aut(F,B_F).
$
For $\sigma\in I(P)$, we have 
$s_\sigma(P)=\{\mbox{id}_F\}$. Therefore
$s_\sigma$ maps $W^0$ into the connected
component of the identity $\Aut^0(F,B_F)$.

We have 
$
s_\sigma(w)\cdot \varphi(x\times_Y \sigma^{-1}w)=
\varphi(x\times_Y w),
$
where $\Phi(x\times_Y w)=(\varphi(x\times_Y w),w)$.
The sections $s_\sigma$ satisfy the identity
$$
s_{\sigma\eta}(w)=s_\sigma(w)\circ
s_\eta(\sigma^{-1}w), \mbox{ for }
\sigma,\eta\in H, w\in W^0.
$$
Therefore they define a $1$-cocycle
$
\xi=\{s_\sigma\}_{\sigma\in H}\in H^1(H, A(W^0)),
$
where $A(W^0)$ is the group of sections of
$\Aut^0(F,B_F)\times W^0$ over $W^0$. 
The $H$-module $\Aut^0(F,B_F)$ is commutative 
by Proposition~\ref{aut}, hence $\xi$ has finite 
order by~\cite{Serre}. After possibly
changing the trivialization, $s_\sigma$ is trivial 
for every $\sigma \in H$.

(5) We have a base change diagram
\[ \xymatrix{
(X,B)|_{Y^0} \ar[d]_f & (F,B_F)\times W'
\ar[d]^{p'}\ar[l] &  (F,B_F)\times W^0    
\ar[l]\ar[d]^p \\
Y^0          & W'\ar[l] & W^0\ar[l]
} \]
where $W'=W^0/H$. The covering $W'\to Y^0$ is
\'etale Galois and $\Gal(W'/Y^0)$ acts on $W'$ 
and $(F,B_F)\times W'$ without fixed points.

The normalization of $Y$ in the field 
$k(W')$ satisfies the required properties.
\end{proof}

\begin{thm} Let $(X,B)$ be a projective log variety 
with Kawamata log terminal singularities such that 
$K+B\sim_\Q 0$. Then 
\begin{itemize}
\item[(i)] the Albanese map $X\to \Alb(X)$
is a surjective morphism, with connected fibers.
\item[(ii)] there exist an \'etale morphism 
$A'\to \Alb(X)$, a projective log variety 
$(F,B_F)$ and an isomorphism over $A'$
$$
(F,B_F)\times A'\isoto (X,B)\times_{\Alb(X)} A'. 
$$
\end{itemize}
\end{thm}

\begin{proof} (i) 
Let 
$
X\stackrel{f}{\to} Y\stackrel{\pi} {\to}A=\Alb(X)
$ 
be the Stein factorization of $\alpha_X$. Since
$f$ is a contraction and $K+B\sim_\Q 0$, 
$f\colon (X,B)\to Y$ is an lc-trivial fibration 
with 
$$
K_Y+\bB_Y+\bM_Y\sim_\Q 0.
$$
Since $\pi$ is finite on its image, any resolution
of $Y$ has non-negative Kodaira dimension 
(\cite{Ueno}, Theorem 10.3). Therefore $K_Y$ is 
$\Q$-linearly equivalent to an effective $\Q$-Weil divisor.
Moreover, $\bB_Y$ is effective since $B$ is effective, 
and $\bM_Y$ is $\Q$-linearly equivalent to an effective 
$\Q$-Weil divisor since $\kappa(\bM)\ge 0$. Therefore
$K_Y\sim_\Q 0$, $\bB_Y=0$ and $\bM_Y\sim_\Q 0$. The
latter implies $\bM\sim_\Q 0$, since $\bM$ is b-nef 
and good.

By Inversion of Ajunction, $Y$ has Kawamata log 
terminal singularities and $K_Y\sim_\Q 0$. The
index one cover $\tilde{Y}\to Y$ has canonical 
singularities and $K_{\tilde{Y}}\sim 0$. Therefore
$\tilde{Y}$, hence $Y$, map onto $A$ by 
\cite{Ueno}, Theorem 10.3.
The finite map $\pi\colon Y\to A$ is \'etale
in codimension one since 
$0\le R=K_Y-\tau^*(K_A)\sim_\Q 0$. Since
$A$ is nonsingular, $\pi$ is \'etale 
everywhere. In particular, $Y$ is an Abelian 
variety. Therefore $\pi$ is an isomorphism, 
by the universality of the Albanese map.

(ii) We have an lc-trivial fibration
$
\alpha_X\colon (X,B)\to A
$ 
with $\bB_A=0$ and $\bM\sim_\Q 0$. By 
Theorem~\ref{str}, there exists a finite 
morphism $\tau\colon A'\to A$, \'etale in 
codimension one, from a normal variety $A'$ 
such that $(X,B)\times_A A'$ is isomorphic to
$(F,B_F)\times A'$ over the generic point of
$A'$. Since $A$ is nonsingular, $\tau$ is 
\'etale everywhere.
Since $\bB_A=0$ and $A'$ is nonsingular, 
we get an isomorphism 
$(X,B)\times_A A'\to (F,B_F)\times A'$
by Proposition~\ref{extend}.
\end{proof}

%%%%%%%%%%%%%%%%%%%%%%%%%%%%%%%%%%%%%%%
%%%%%%%%%%%%%%%%%%%%%%%%%%%%%%%%%%%%%%%

\end{document}